\documentclass[reqno,12pt]{amsart}
\usepackage{latexsym, amsfonts,mathrsfs, amsmath, amssymb, amscd, epsfig}
\usepackage{mathrsfs}

\newtheorem{thm}{Theorem}
\newtheorem{defn}[thm]{Definition}
\newtheorem{prop}[thm]{Proposition}
\newtheorem{lem}[thm]{Lemma}
\newtheorem{cor}[thm]{Corollary}
\newtheorem{rmk}{Remark}

\newenvironment{pf}{{\noindent \it \bf Proof:}}{{\hfill$\Box$}\\}

\def\beq#1\eeq{\begin{equation}#1\end{equation}}
\def\balign #1 #2 \ealign{\begin{aligned} #1 #2  \end{aligned} }

\def\ep{\varepsilon}
\def\s{\sigma}

\def\bX{\mathbf{X}}
\def\mcE{\mathcal{E}}
\def\hmcE{\hat{\mathcal{E}}}
\def\tmcE{\tilde{\mathcal{E}}}

\def\brhop{\bar{\rho}_+}

\def\bMp{\bar{M}_+}
\def\bup{\bar{u}_+}
\def\bum{\bar{u}_-}

\def\tx{\tilde{x}}

\def\bPsi{\bar{\Psi}}

\def\hmD{\hat{\mathscr{D}}}

\def\mL{\mathscr{L}}

\def\mA{\mathscr{A}}
\def\mC{\mathscr{C}}

\def\mF{\mathscr{F}}

\def\mD{\mathscr{D}}

\def\tx{\tilde{x}}

\def\lbn{\||}
\def\rbn{|\|}
\def\ttn{\tilde{\|}}
\def\tn{\tilde{|}}

\def\mBV{\mathbb{V}}
\def\mBR{\mathbb{R}}

\begin{document}

\title[ ]{Global Stability of Steady Transonic Euler Shocks
 in  Quasi-One-Dimensional Nozzles}

\author[ ]{Jeffrey Rauch}
\author[ ]{Chunjing  Xie}
\address{Department of Mathematics, University of Michigan, 530 Church
Street, Ann Arbor, MI 48109, USA} \email{rauch@umich.edu}
\email{cjxie@umich.edu}
\author[ ]{Zhouping Xin}
\address{The Institute of Mathematical Sciences and department of
mathematics, The Chinese University of Hong Kong, Hong Kong}
\email{zpxin@ims.cuhk.edu.hk}


\begin{abstract}
We prove global in time dynamical stability of
steady transonic shock solutions in  divergent quasi-one-dimensional nozzles.
We assume
neither the smallness of the relative slope of the nozzle nor the weakness of the shock.
 Key ingredients of the proof are  an exponentially decaying
energy estimate for a
linearized problem together with  methods from \cite{LRXX}.

\end{abstract}

\maketitle

\section{ Introduction and Main Results}

Compressible isentropic Euler flows in quasi-one-dimensional
nozzles are governed by
\begin{equation}\label{Q1DNozzle}
\begin{aligned}
& \rho_t+(\rho u)_x=-\frac{a^{\prime}(x)}{a(x)}\rho u, \\
& (\rho u)_t+ (\rho u^2 +p(\rho))_x=-\frac{a^{\prime}(x)}{a(x)}\rho u^2,
\end{aligned}
\,.
\end{equation}
where $\rho, u$ and $p$ denote  the density, velocity
and pressure, and, $a(x)$ is the cross-sectional area of the nozzle. The typical examples are 1D,  2D rotationally symmetric, and  3D spherical symmetric compressible Euler equations, where $a(x)$ are $1$,  $x$  and $x^2$, respectively.

We assume that $p$ satisfies
\begin{equation*}
\forall \rho>0, \qquad  p(\rho)>0\,,\quad p'(\rho)>0\,,
\quad \text{and}\,,\quad p''(\rho)\geq 0\,.
\end{equation*}
Examples are  $p(\rho):=A\rho^{\gamma}$ with $\gamma\geq 1$,
in particular the  {\it polytropic} ($\gamma>1$) and {\it isothermal} ($\gamma=1$) gases. The {\it sound speed} is $c(\rho):=\sqrt{p'(\rho)}$. A state $\rho,u$
at which the flow speed is larger than the sound speed,
$|u|>\sqrt{p'(\rho)}$,  is  called {\it supersonic}. The opposite case, $|u|<\sqrt{p'(\rho)}$,  is {\it subsonic}.

For steady states,   (\ref{Q1DNozzle}) becomes the  system of ordinary differential
equations
\begin{equation}\label{SteadyQ1D}
\begin{aligned}
&(\rho u)_x=-\frac{a^{\prime}(x)}{a(x)}\rho u, \\
&(\rho u^2 +p(\rho))_x=-\frac{a^{\prime}(x)}{a(x)}\rho u^2\,.
\end{aligned}
\,,
\end{equation}
studied intensively in  \cite{EmbidGM, LiuARMA}.  Some
properties needed later are the following.

The first equation in (\ref{SteadyQ1D}), $(a\rho u)_x=0$, holds if and only if
\begin{equation}\label{massflux}
\exists M\in {\mathbb R},\qquad a\rho u\ =\ M\,.
\end{equation}
 After  a possible reflection we may assume without
 loss of generality that $M>0$. The flow is from left to right.
Using   (\ref{massflux}) to eliminate $u$ from the second equation in \eqref{SteadyQ1D} yields
 a scalar equation for $\rho$,
\begin{equation}\label{Steadymomeq}
\left(\frac{M^2}{a\rho}\right)_x\ +\ a p(\rho)_x\ =\ 0.
\end{equation}

\begin{defn}
A {\it  steady transonic shock solution} is a piecewise smooth solution of (\ref{Steadymomeq})
with two smooth solutions separated by  a shock connecting a supersonic state on the left
to a subsonic state on the right.
\end{defn}

The steady transonic shock solutions of (\ref{SteadyQ1D}) on $[l, L]$ are the
piecewise smooth solutions
\begin{equation}
\label{STshock}
\begin{aligned}
(\bar\rho,\bar{u})=\left\{
\begin{array}{ll}
(\bar\rho_-,\bar{M}/(a\bar{\rho}_-)),\quad\text{if}\,\,l<x<x_0,\\
(\bar\rho_+,\bar{M}/(a\bar{\rho}_+)),\quad\text{if}\,\,x_0<x<L,\\
\end{array}
\right.
\end{aligned}
\end{equation}
satisfying the Rankine-Hugoniot condition
\begin{equation}\label{bRH}
\left(p(\bar\rho_-)+\frac{\bar M}{(a\bar\rho_-)^2}\right)(x_0-)
\ =\
 \left(p(\bar\rho_+)+\frac{\bar M}{(a\bar\rho_+)^2}\right)(x_0+),
\end{equation}
with $\bar\rho_-$ is supersonic ($p'(\bar\rho_-)< \frac{\bar M^2}{(a\bar\rho_-)^2}$) and $\bar\rho_+$ is subsonic ($p'(\bar\rho_+)< \frac{\bar M^2}{(a\bar\rho_+)^2}$). This
implies  the entropy condition
\begin{equation}
\bar{\rho}_-(x_0)\ <\  \bar{\rho}_+(x_0)\,.
\end{equation}

Steady transonic shock solutions were constructed in \cite{EmbidGM, LiuARMA}
satisfying the boundary conditions
\begin{equation}\label{SimpleBC}
 (\rho, u)(l)\ =\ (\rho_l, u_l)\ =\ (\rho_l, \frac{\bar M}{a(l)\rho_l}),\quad
 \text{ and, }
 \quad
 \rho(L)\ =\ \rho_r,
\end{equation}
with  $(\rho_l, u_l)$ and $(\rho_r,u_r)$)  supersonic and subsonic respectively.
There have been many studies of
the stability of steady transonic shocks  for  (\ref{Q1DNozzle}) (see
\cite{LiuTP}).
In \cite{LiuTP}
 a wave front tracking variant of Glimm's
scheme was used to prove that
when $|a'(x) / a(x) | $  is small,
a weak transonic shock  is dynamically stable if $a'(x_0)>0$ and dynamically unstable if $a'(x_0)<0$.
Studies on the solutions of general
hyperbolic conservation laws with source terms include \cite{LiuARMA,
LiuJMP, Lien, Ha, HaYang} and references therein. For piecewise
smooth initial data, in the case that $a(x)=x$ or $a(x)=x^2$,
Xin and Yin \cite{XinYinJDE}
proved dynamical stability of weak shocks in nozzles
 with large $l$ (thus $ | a'(x)/a(x | $   small) in the class of piecewise
 smooth transonic shock solutions.
The main  result of this paper  proves
global stability in time without these restrictions.

Let $(\bar\rho, \bar u)$ be a steady transonic shock solution of  the form (\ref{STshock}) satisfying the boundary conditions (\ref{SimpleBC}). Suppose that  the solution stays away from vacuum,
\begin{equation}\label{vacuum}
\inf_{x\in [l, L]}  \   \bar \rho(x)\ >\
0\,.
\end{equation}
Solving an initial value problem for
\eqref{SteadyQ1D}
we can extend $(\bar\rho_-, \bar u_-)$ to be a smooth supersonic solution to (\ref{SteadyQ1D}) on $[l, x_0+\delta]$ for some $\delta>0$, that
 coincides with $(\bar\rho_-, \bar u_-)$ on $[l, x_0]$.
The notation $(\bar\rho_-, \bar u_-)$ is used for
the extended solution as well. Similarly, denote by $(\brhop, \bup)$
 a subsonic solution of (\ref{SteadyQ1D}) on $[ x_0-\delta, L]$ for some $\delta>0$,
 that coincides with $(\brhop, \bup)$  on $[x_0, L]$.

Consider the mixed  initial boundary value problem defined by
(\ref{Q1DNozzle})  with initial data,
\begin{equation}
\label{initial}
(\rho, u)(0, x)\ =\
 (\rho_0,u_0)(x),
 \end{equation}
and boundary conditions
\begin{equation}
\label{bc2}
(\rho, u)(t, l)\ =\
(\rho_l,\frac{\bar M}{a(l)\rho_l}),\qquad
 \rho(t, L)\ =\ \rho_r\,,
 \end{equation}
with $\rho_l$ and $\rho_r$ from (\ref{SimpleBC}).

Assume that the initial data satisfy
\begin{equation}\label{in1}
(\rho_0,u_0)(x)=\begin{cases}
&(\rho_{0-},u_{0-})(x),\quad\text{if} \ \  l<x<\tilde{x}_0,\\
&(\rho_{0+},u_{0+})(x),\quad\text{if}  \ \ \tilde{x}_0<x<L,
\end{cases}
\end{equation}
These data are  small perturbations of $(\bar\rho, \bar u)$ in the sense that
\begin{equation}
\label{AssID}
\begin{aligned}
|x_0-\tilde{x}_0|  &\ +\  \big\|
(\rho_{0+},u_{0+})\  -\ (\bar\rho_+,\bar u_+)
\big\|_{H^{k+2}([\check{x}_0, L])} \\
&+ \big\|
(\rho_{0-},u_{0-}) -(\bar\rho_-,\bar u_-)\big\|_{H^{k+2}([l,\hat{x}_0])}<\varepsilon,
\end{aligned}
\end{equation}
for some small $\varepsilon>0$, and some integer $k\ge 15$, where $\check{x}_0=\min\{x_0, \tilde{x}_0\}$ and $\hat{x}_0=\max\{x_0, \tilde{x}_0\}$. Moreover, $(\rho_0, u_0)$ is assumed to satisfy the Rankine-Hogoniot condition.
That is,  at the shock location  $x=\tilde{x}_0$,
\begin{equation}
\label{rh2}
 \big((p(\rho_{0+})+\rho_{0+}u_{0+}^2   - (p(\rho_{0-})+\rho_{0-}u_{0-}^2)\big)
\,(\rho_{0+}-\rho_{0-})
\ =\ (\rho_{0+}u_{0+}-\rho_{0-}u_{0-})^2.
\end{equation}
Our dynamical stability theorem  is the following.

\begin{thm}
\label{Thmdystability}
Let $(\bar\rho,\bar{u})$ be a steady transonic shock solution to system (\ref{SteadyQ1D}) satisfying (\ref{STshock}),  (\ref{bRH}), (\ref{SimpleBC}), and (\ref{vacuum}).   The key hypothesis is
that the nozzle is widening at the shock location,
\begin{equation}
\label{stabilitycondition}
a'(x_0)\ >\ 0.
\end{equation}
There exists an $\ep_0>0$ so that for any $0<\ep\leq \ep_0$,
if the initial data $(\rho_0, u_0)$ satisfy  (\ref{AssID}), (\ref{rh2})
and the  compatibility conditions at $x=l$, $x=\tilde{x}_0$ and $x=L$ or order
$(k+2)$,
then the initial boundary value problem (\ref{Q1DNozzle}), (\ref{initial}) and (\ref{bc2})
has  a unique piecewise smooth solution $(\rho, u)(x, t)$ for $(x, t)\in [l, L]\times [0, \infty)$ containing a single transonic shock $x=s(t)$ with $s(0)=\tilde{x}_0$  and $l<s(t)<L$ satisfying the
pair of Rankin-Hougoniot conditions for $t\ge 0$,
$$
(p(\rho)+\rho u^2)(t,s(t+))-(p(\rho)+\rho u^2)(t, s(t-))\
 =\ \big(\rho u(t, s(t+))-\rho u(t, s(t-))\big)\dot s(t),
$$
$$
\rho u(t, s(t+))-\rho u(t, s(t-))
\ =\ \big(\rho (t, s(t+))-\rho(t, s(t-)\big))\dot{s}(t),
$$
and   the Lax geometric shock conditions,
$$
(u-\sqrt {p'(\rho)})(t, s(t-))
\ >\
\dot s(t)\ >\
(u-\sqrt {p'(\rho)})(t, s(t+)),
\qquad
(u+\sqrt {p'(\rho)})(t, s(t+))
\ >\
\dot s(t).
$$
Denote
\begin{equation*}
(\rho,u)=\left\{
\begin{array}{ll}
(\rho_-,u_-),\quad\text{if}\quad l<x<s(t),\\
(\rho_+,u_+),\quad\text{if}\quad s(t)<x<L\,.
\end{array}
\right.
\end{equation*}
There is a $T_0>0$ so that for $t>T_0$,
\begin{equation*}
\begin{aligned}
(\rho_-,u_-)(t, x)=(\bar\rho_-,\bar u_-)(x), \quad \text{for}\,\, l\le x<s(t)\,.
\end{aligned}
\end{equation*}
The solution approaches the steady transonic flow at an exponential rate, that is,
there exist positive  constants $C>0$ and $\lambda>0$
independent of the solution so that
for $t>0$
\begin{equation*}
\begin{aligned}
\|(\rho_+,u_+)(\cdot,t)-(\bar\rho_+,\bar u_+)(\cdot)\|_{W^{k-7,\infty}(s(t), L)} +\sum_{m=0}^{k-6}|\partial^{m}_t(s(t)-x_0)|
\ \leq\
 C\,\varepsilon \,e^{-\lambda t},
\end{aligned}
\end{equation*}
 where  $(\bar{\rho}_{\pm}, \bar{u}_{\pm})$ is the unperturbed
solution.
\end{thm}

\begin{rmk}
 The results in  Theorem \ref{Thmdystability}, are also true if we
impose small perturbations in  the boundary conditions
(\ref{SimpleBC}).
\end{rmk}

\begin{rmk}
The regularity assumption in (\ref{AssID}) is not optimal. Adapting the paradifferential
methods from
\cite{Metivier}, one can decrease the regularity required  in
(\ref{AssID}).
\end{rmk}

\begin{rmk}
Compatibility conditions for the initial boundary value problems for hyperbolic equations
are discussed in detail in \cite{RauchM, Majda, Metivier}.
\end{rmk}

\begin{rmk}
We require neither the  smallness of $ | a^\prime / a | $,  largeness of $l$,
nor the  weakness of the shock strength as  in \cite{XinYinJDE}.  Similar results for
steady  transonic shock solutions for the Euler-Poisson equations were
proved in \cite{LRXX}. The strategy of the  proof
of  Theorem 1 is  inspired by  \cite{LRXX}.
\end{rmk}

\begin{rmk}
If the condition (\ref{stabilitycondition}) is violated, i.e. $a^\prime(x_0)<0$, then the transonic shock solution is unstable, see \cite{LiuTP, XinYinJDE}.
\end{rmk}

The paper is organized as follows.
Section 2 transforms the problem to
 free boundary  problem for a second order scalar
 hyperbolic partial differential equation.
 The weighted energy equality for the associated linearized problem is derived in \S 3.  This yields Theorem 1 with the aid of
 ideas from  \cite{LRXX}.


\section{Transformation of the Problem}
Let $(\bar{\rho},\bar{u})$ be a steady transonic shock solution   of the form (\ref{SteadyQ1D}) satisfying (\ref{STshock}), (\ref{bRH}), (\ref{SimpleBC}),  and (\ref{vacuum}).
If  the initial data $(\rho_0,u_0)$  satisfies (\ref{AssID}) and the compatibility condition
of order $k+2$.   It follows from the argument in \cite{LY}
that there exists a piecewise smooth solution  of the
equations on $[0, \bar T]$ for some $\bar T>0$, which can be written as
\begin{equation}\label{ustsolution}
(\rho,u)=\left\{
\begin{array}{ll}
(\rho_-,u_-),\,\,\text{if}\,\,l<x<s(t),\\
(\rho_+,u_+),\,\,\text{if}\,\,s(t)<x<L.
\end{array}
\right.
\end{equation}
Since the flow is supersonic on the left it follows that
there is a $T_0>0$ so that, when $t>T_0$, $(\rho_-, u_-)$ depends
 only on the boundary conditions at $x=l$. Moreover, when $\varepsilon$ is small, by the standard lifespan argument, we have $T_0< \bar T$
(cf. \cite{LY}). Therefore,
\begin{equation}\label{ste}
(\rho_-, u_-)=(\bar{\rho}_-, \bar{u}_-)\,\, \text{for}\,\, t>T_0.
\end{equation}
It suffices to study the solution for $T_0>0$ so
 without loss of generality we suppose that  $T_0=0$.
We want to extend the local-in-time solution to all $t>0$. In view of (\ref{ste}),
it suffices to obtain uniform estimates in the region $x>s(t)  ,\, t>0$.
We
formulate an initial boundary value problem in this region.  First, the Rankine-Hugoniot conditions for (\ref{ustsolution})  read
\begin{equation}
\label{USRH}
[\rho u]\ =\ [\rho]s'(t),\qquad
[\rho u^2+p]\ =\ [\rho u]s'(t),
\end{equation}
using standard notation for the jump
$[f] := ,f(s(t+),t)-f(s(t-),t)$, so
\begin{equation*}
[p+\rho u^2] \ [\rho]\ =\
[\rho u]^2.
\end{equation*}
Denote $M\ :=\ a\rho u$. Then
\begin{equation*}
\begin{aligned}
&\left(p(\rho_+)(t,s(t))+\frac{M_+^2(t,s(t))}{a^2(s(t))\rho_+(s(t),t)}
-p(\bar\rho_+)(s(t))
-\frac{\bar
M_+^2}{a^2(s(t))\bar\rho_+}(s(t))\right.\\
&\left.+\,p(\bar\rho_+)(s(t))+\frac{\bar M_+^2}{a^2\bar\rho_+}(s(t))
-p(\bar\rho_-)(s(t))-\frac{\bar
M_-^2}{a^2 \bar\rho_-}(s(t))\right)\ \big(\rho_+(t,s(t))-\rho_-(s(t))\big)\\
&\  = \left(\frac{M_+(t,s(t))-\bar M_+(s(t))}{a(s(t))}\right)^2.
\end{aligned}
\end{equation*}
The implicit function theorem and the momentum equation in (\ref{SteadyQ1D})
imply,
\begin{equation}\label{DiffJshock}
(M_+-\bMp)(t,s(t))=\mA_1((\rho_+-\brhop)(t,s(t)), s(t)-x_0)
\end{equation}
where $\mA_1$ regarded as a function of two variables satisfies $\mA_1(0,0)=0$ and
\begin{equation*}
\begin{aligned}
 \frac{\partial \mA_1}{\partial
(\rho_+-\brhop)} \Bigg|_{(0,0)}=-\frac{a(p'(\bar\rho_+)-\bup^2)}{2\bup}(x_0),\qquad
\frac{\partial \mA_1}{\partial (s-x_0)}\Bigg|_{(0,0)}=-
\frac{a'\bum(\bar\rho_+-\bar\rho_-)}{2}(x_0).
\end{aligned}
\end{equation*}
Note that $\bar M_+=\bar{M}_-=\bar M$. And so  substituting
(\ref{DiffJshock}) into the first equation in (\ref{USRH})
yields,
\begin{equation}\label{shockspeedeq}
s'(t)
\ =\
\mA_2(\rho_+-\brhop,s(t)-x_0)
\end{equation}
where $\mA_2$ satisfies $\mA_2(0,0)=0$ and
\begin{equation*}
\frac{\partial \mA_2}{\partial(\rho_+-\brhop)}\Bigg|_{(0,0)}=-\frac{p'(\bar\rho_+)-\bup^2}{2\bar
u_+(\bar\rho_+-\bar\rho_-)}(x_0),
\qquad
\frac{\partial \mA_2}{\partial(s(t)-x_0)}\Bigg|_{(0,0)}=-\frac{a'\bum}{2a}(x_0).
\end{equation*}

Define
\begin{equation*}
\begin{aligned}
\psi_-(t, x)\ :=\
-\bar{M}t +\int_l^{x} a(y) \bar\rho_-(y)dy
\end{aligned}
\end{equation*}
and
\begin{equation*}
\begin{aligned}
\bar{\psi}_+(t, x)\ :=\ {\psi}_-(t, x_0)+\int_{x_0}^x a(y)\bar{\rho}_+(y)dy,\,\, \psi_+(t, x)=\psi_-(t, s(t))+ \int_{s(t)}^x a(y)\rho_+(y, t)dy.
\end{aligned}
\end{equation*}
Set
\begin{equation*}
\begin{aligned}
\Psi
\ :=\
\psi_+(x,t)-\bar \psi_+(t, x).
\end{aligned}
\end{equation*}
Then
\begin{equation*}
\Psi_t=\bar M-M_+,\qquad
\Psi_x=a\rho_+-a\brhop.
\end{equation*}
The second equation in (\ref{Q1DNozzle}) is equivalent to
\begin{equation*}
\begin{aligned}
-(\bar{\psi}_++\Psi)_{tt}+\partial_x\left(\frac{(\partial_t\bar{\psi}_+
+\partial_t\Psi)^2}
{\partial_x\bar{\psi}_++\partial_x\Psi}\right)+a(x)\partial_xp\left(
\frac{\partial_x\bar{\psi}_++\partial_x\Psi}{a(x)}\right)=0.
\end{aligned}
\end{equation*}
Write this equation as follows. Setting $\xi=(\xi_0, \xi_1)=(t, x)$,
\begin{equation}\label{EQY}
\sum_{ij}\bar{a}_{ij}(x,\Psi_t,\Psi_x)\partial_{ij}\Psi+\sum_i\bar{b}_i(x,\Psi_t,\Psi_x)\partial_i
\Psi =0,
\end{equation}
where $a_{ij}$, $b_i$ and $c$ are smooth with respect to each variable, and satisfy
(using the Einstein summation convention),
\begin{equation*}
\begin{aligned}
\mL_0 \Psi &= \bar a_{ij}(x,0,0)\partial_{ij}\Psi+ \bar b_i(x,0,0)\partial_i \Psi+c \bar g(x,0,0)\Psi\\
& = \partial_{tt}\Psi+2\bup
\partial_{tx}\Psi+(\bup^2-p'(\brhop))\partial_{xx}\Psi +2\partial_x
\bup \partial_t \Psi +B\partial_x\Psi
\end{aligned}
\end{equation*}
with
\begin{equation*}
B
\ :=\
\frac{a'}{a}p'(\brhop)-\partial_xp'(\brhop)+\partial_x\bup^2.
\end{equation*}

In terms of $\Psi$, the equations (\ref{DiffJshock}) and
(\ref{shockspeedeq}) can be written as
\begin{equation}\label{1stUSRH}
\Psi_t(t, s(t))=-\mA_1(\frac{\Psi_x(t, s(t))}{a(s(t))}, s(t)-x_0),
\end{equation}
and
\begin{equation}\label{shockspeedY}
\begin{aligned}
s'=\mA_2(\frac{\Psi_x(t, s(t))}{a(s(t))},s(t)-x_0),
\end{aligned}
\end{equation}
respectively.
A direct computation yields
\begin{equation*}
\begin{aligned}
\Psi(t, s(t))&=\psi_+(t, s(t))-\bar \psi_+(t, s(t))
= \psi_-(t, s(t))-\bar \psi_+(t, s(t))\\
&= \psi_-(t,x_0)+\psi_-(t,s(t))- \psi_-(t,x_0)-\bar \psi_+(t,x_0)-\bar
\psi_+(t,s(t))+\bar \psi_+(t,x_0)\\
&=\partial_x\psi_-(t,x_0)\cdot(s(t)-x_0)-\partial_x\bar
\psi_+(t,x_0)\cdot(s(t)-x_0)+R_1,
\end{aligned}
\end{equation*}
where $R_1$ is quadratic in $s(t)-x_0$. This implies
\begin{equation}\label{shockY}
s(t)-x_0=\mA_3(\Psi(s(t),t)),
\end{equation}
where $\mA_3$ satisfies $\mA_3(0)=0$ and
\begin{equation*}
\frac{\partial\mA_3}{\partial \Psi}\Bigg|_{\Psi=0}
\ =\ \frac{1}{a(\bar{\rho}_- -\brhop)}(x_0)\,.
\end{equation*}
It follows from (\ref{1stUSRH}), (\ref{shockspeedY}) and (\ref{shockY}) that
\begin{equation}\label{BCYshock}
\begin{aligned}
\partial_t \Psi=\mA_4(\Psi_x, \Psi),\,\,\,\,\,\,\text{at}\,\,x=s(t),
\end{aligned}
\end{equation}
where
\begin{equation*}
\mA_4(0,0)=0,
\qquad
 \frac{\partial \mA_4}{\partial
\Psi_x}=\frac{c^2(\bar\rho_+)-\bar u_+^2}{2\bar
u_+}(x_0),\qquad
\frac{\partial \mA_4}{\partial
\Psi}=-\frac{a'\bar{u}_-}{2a}(x_0)\,.
\end{equation*}
On the right boundary, $x=L$, $\Psi$ satisfies
\begin{equation}
\label{BCYR}
\partial_x \Psi=0 ,
\qquad
\text{at}\
\quad
x=L \,.
\end{equation}
It suffices to derive uniform estimates for $\Psi$ and $s$ that
satisfy (\ref{EQY}), (\ref{shockY}), (\ref{BCYshock}) and
(\ref{BCYR}).

To this end,
first transform this free boundary value problem into a fixed boundary
value problem.
Set
\begin{equation*}
\tilde{t}:=t,\qquad
 \tilde{x}:=(L-x_0)\frac{x-s(t)}{L-s(t)}+x_0\,,
 \qquad
 \sigma(\tilde{t}):=s(t)-x_0,
\end{equation*}
and
\begin{equation}\label{q1q2}
q_1(\tilde x, \s):=\frac{L-\tilde x}{L-x_0\ -\ \sigma (\tilde t)}\,,
\qquad
q_2(\s)=\frac{L-x_0}{L-x_0-\sigma (\tilde t)}\,.
\end{equation}
Then
$$
\frac{\partial}{\partial t}=\frac{\partial}{\partial \tilde t}
\ -\ \sigma '(\tilde t)q_1\frac{\partial}{\partial{\tilde x}},\qquad \frac{\partial}{\partial x}=q_2\frac{\partial}{\partial{\tilde x}},$$
\begin{align*}
\frac{\partial^2}{\partial t^2}&=\frac{\partial^2}{\partial \tilde t^2}+(q_1 \sigma' (\tilde t))^2\frac{\partial^2}{\partial \tilde x^2}-2q_1\sigma'(\tilde t)\frac{\partial^2}{\partial \tilde x\partial \tilde t}-q_1\left(\sigma''(\tilde t)+2\frac{( \sigma' (\tilde t))^2}{L-x_0-\sigma (\tilde t)}\right)\frac{\partial}{\partial \tilde x},\end{align*}
$$\frac{\partial^2}{\partial x\partial t}=q_2\left (\frac{\partial^2}{\partial \tilde x\partial \tilde t}+\frac{\sigma'(\tilde t)}{L-x_0-\sigma (\tilde t)}\frac{\partial}{\partial{\tilde x}}-q_1\sigma'(t)\frac{\partial^2}{\partial \tilde x^2}\right),
\qquad
\frac{\partial^2}{\partial  x^2}=q_2^2\frac{\partial^2}{\partial \tilde x^2}.
$$
So (\ref{EQY}) becomes
\begin{equation}
\label{eqinnew}
\begin{aligned}
\partial_{\tilde t\tilde t}\Psi\   +\
 &  q_2\partial_{\tilde x}\left( \frac{(-\bar M+\partial_{\tilde{t}} \Psi -\s'(\tilde t)
q_2 \partial_{\tilde{x}} \Psi)^2}{a\brhop
\ +\
q_2\partial_{\tilde{x}}\Psi}\right)+a(x) q_2 \partial_{\tilde{x}}
p\left(\frac{a\brhop+q_2\Psi_{\tilde{x}}}{a(x)}\right) \notag
\\
&
-\ 2\sigma'(t) q_1 \partial_{\tilde x\tilde
t}\Psi
\ +\ (q_1\sigma'(t))^2\partial_{\tilde x \tilde x}\Psi
\ -\
2\frac{(\sigma'(\tilde t))^2}{L-x_0-\sigma(\tilde t)}q_1 \partial_{\tilde x}\Psi
\  =\
\sigma''(\tilde t)q_1 \partial_{\tilde x}\Psi.
\end{aligned}
\end{equation}

Equation \eqref{shockY}
takes the form
\begin{equation*}
\s = \mA_3(\Psi(t, \tx =x_0)),
\end{equation*}
and the equation for the shock front,  (\ref{shockspeedY}),  becomes
\begin{equation*}
\frac{d\sigma}{d\tilde{t}}=\mA_2\left(\frac{q_2(\s)\Psi_{\tx}}{a},\ \sigma(\tilde{t})\right).
\end{equation*}
Using (\ref{shockY}) to represent the quadratic terms for
$\sigma$ in terms of $\Psi$, we have, at $\tilde x=x_0$,
\begin{equation}\label{shockseq}
\frac{d\sigma}{d\tilde{t}}+\frac{a'\bar
u_-}{2a}(x_0)\sigma=\mC_2(\Psi_{\tilde{x}}, \Psi),
\end{equation}
where $\mC_2$ satisfies
\begin{equation*}
\left| \mC_2(\Psi_{\tilde{x}}, \Psi) +
\frac{c^2(\brhop)-\bup^2}{2(\brhop-\bar{\rho}_-)\bup a}
(x_0)\Psi_{\tilde{x}}    \right| \leq C(\Psi_{\tilde{x}}^2+\Psi^2).
\end{equation*}
It follows from (\ref{shockY}) and (\ref{shockseq}) that
one can represent $\s$ and $\s'$ in terms of $\Psi$ and its
derivatives at $\tilde{x}=x_0$. Thus, after manipulating
(\ref{BCYshock}) with (\ref{shockY}) and (\ref{shockseq}), one finds
\begin{equation*}\label{defC1}
\Psi_{\tilde {t}}
\ =\ \mC_1(\Psi_{\tx}, \Psi), \qquad {~\rm at~} \quad \tx=x_0.
\end{equation*}
By the implicit function theorem this is equivalent to
\begin{equation*}
\Psi_{\tx}=\mC_3(\Psi_{\tilde{t}}, \Psi), \qquad \text{at} \quad
\tx=x_0,
\end{equation*}
where $\mC_3$ satisfies
\begin{equation*}
\left|\mC_3(\Psi_{\tilde{t}}, \Psi)-
\frac{2\bup}{c^2(\brhop)-\bup^2}(x_0)\Psi_{\tilde{t}}-\frac{a'\bar
u_- \bup}{(c^2(\brhop)-\bup^2) a}(x_0)\Psi\right|
\ \leq\
C(\Psi_{\tilde{t}}^2+\Psi^2).
\end{equation*}

In the following,  drop the $\ \tilde{  }\ $ in $\tilde x$ and $\tilde t$ for
ease of reading.

In summary, the problem has been transformed to  the following compact form
\begin{equation}
\label{nonlinearpb}
\left\{
\begin{aligned}
&\mL(x,{\Psi},\sigma){\Psi} \ =\ \sigma''(t)q_1 \partial_{ x}\Psi, \quad
(t,x)
\in [0,\infty)\times [x_0, L],\\
&\partial_{x}{\Psi}\ =\
d_1(\Psi_t, \Psi)\Psi_t+e_1(\Psi_t, \Psi) \Psi,\quad \text{at} \,\, x=x_0,\\
&\partial_{x} \Psi\ =\ 0,\quad
\text{at}\,\, x=L,\\
&\sigma({t}) \ =\ \mA_3(\Psi(t, x_0)), \quad \s(0)=\s_0,
\end{aligned}
\right.
\end{equation}
where,
\begin{equation*}
\begin{aligned}
\mL(x,\Psi, \sigma)\Phi
\  :=\
&\sum_{i,j=0}^1a_{ij}(x, \nabla
\Psi,\sigma, \sigma'
)\partial_{ij}\Phi+\sum_{i=0}^1b_i(x, \nabla \Psi,\sigma,
\sigma')
\partial_i \Phi,
\end{aligned}
\end{equation*}
with
\begin{equation*}
\begin{aligned}
&d_1(\Psi_t, \Psi\ :=\  \int_0^1\frac{\partial \mC_3}{ \partial \Psi_t}
( \theta \Psi_t,\theta \Psi)d\theta,
\qquad
e_1( \Psi_t, \Psi)
\ :=\
\int_0^1\frac{\partial \mC_3}{\partial \Psi} ( \theta
\Psi_{t},\theta \Psi)d\theta.
\end{aligned}
\end{equation*}
Furthermore, one has
\begin{equation*}
\begin{aligned}
&\mL(x,0,0)\Phi=\mL_0\Phi
\end{aligned}
\end{equation*}
and
\begin{equation*}
\begin{aligned}
&a_{00}(x,  \nabla \Psi, \sigma, \sigma')=1,\qquad
a_{01}(x,0,0,0,0)=a_{10}(x,0,0,0,0)=\bar u_+,\,\,\\
&a_{11}(x,0,0,0,0)=-(p'(\brhop)-\bup^2),\\
&b_0(x,0,0,0, 0)=2\partial_x\bup, \qquad b_1(x,0,0,0, 0)=B,\\
&d_1(0, 0)=\frac{2\bar
u_+}{c^2(\bar\rho_+)-\bar u_+^2}(x_0),\qquad
e_1(0, 0)= \frac{a'(x_0)\bup \bum}{(c^2(\brhop)-\bup^2)a}(x_0).\end{aligned}
\end{equation*}

\section{Proof of the Main Theorem}
The key first step in the proof of Theorem \ref{Thmdystability}
is the  exponential decay of solutions of a linearized problem.
It yields in non trivial manner {\it  a priori} estimates for the nonlinear problem
which lead to the global  stability of the steady shock.

\subsection{Linear Estimate.}
The first step is the following energy identity for the linearization
at the steady shock.

\begin{lem}\label{keylem}
Assume that $a$ satisfies (\ref{stabilitycondition}). Let $\Psi$ be
a smooth solution of the  linearized problem
\begin{equation}\label{lzeropb}
\left\{
\begin{aligned}
&\mL_0\Psi\ =\ 0,\qquad x_0<x<L, \quad  t>0,
\\
&\partial_x\Psi
\ =\
 \frac{2\bar
u_+}{c^2(\bar\rho_+)-\bar u_+^2}(x_0) \partial_t\Psi
\ +\
 \frac{a'(x_0)\bup \bum}{(c^2(\brhop)-\bup^2)a}(x_0)\Psi,\quad {\rm at}\ \ x=x_0,\\
&\partial_x\Psi=0 ,
\quad
{\rm at} \quad x=L,
\\
&\Psi(0,x)\ =\
h_1(x), \quad
 \Psi_t(0,x)\ =\ h_2(x),\quad  x_0<x<L.
\end{aligned}
\right.
\end{equation}
Then the following dissipation identity holds,
\begin{equation}\label{energyED}
E(\Psi, t)\ =\
 E(\Psi, 0)-D(\Psi, t) ,
\end{equation}
where $E$ and $D$ are defined as,
\begin{equation*}
\begin{aligned}
&E(\Psi, t)\ :=\
\left(\frac{a'\bup^2\bum\Psi^2}{a}\right) (t, x_0)
+\int_{x_0}^L\bup
\left\{(\partial_t \Psi)^2+ (p'(\bar{\rho}_+)-\bup^2)(\partial_x
\Psi)^2 \right\}(t,x) \ dx\,, \\
&D(\Psi, t)\ :=\ 2\int_0^t \bup^2(x_0)(\partial_t\Psi)^2(\tau, x_0)
+\bup^2(L)(\partial_t\Psi)^2(\tau, L)\ d\tau.
\end{aligned}
\end{equation*}
\end{lem}

\begin{pf}  Multiplying the  first equation in (\ref{lzeropb})
by  $\bar{u}_+(x)\partial_t\Psi$
and integrating by parts yields
\begin{equation}
\label{energy1}
\begin{aligned}
\int_0^t \int_{x_0}^L \mL_0 \Psi\cdot \bup \partial_t\Psi \ dxd\tau
&\ =\ \int_{x_0}^L  \frac{\bup(\partial_t\Psi)^2
+\bup (c^2(\brhop-\bup^2) (\partial_x\Psi)^2}{2}(t, x)\ dx\\
&+ \int_0^t \int_{x_0}^L [B \bup +\partial_x(\bup(c^2(\brhop)-\bup^2))]
\partial_t\Psi\partial_x\Psi \ dx d\tau\\
& + \int_0^t \Big(\bup^2(\partial_t\Psi)^2 + \bup(\bup^2-c^2(\brhop))\partial_x \Psi \partial_t\Psi\Big)\Big|_{x=x_0}^{x=L} \ d\tau\\
& - \int_{x_0}^L  \frac{\bup(\partial_t\Psi)^2
+\bup (c^2(\brhop-\bup^2) (\partial_x\Psi)^2}{2}(0, x)\ dx\,.
\end{aligned}
\end{equation}
Note that
\begin{equation}
\label{term1}
\begin{aligned}
B\bup\  +\
&\partial_x(\bup(x)(c^2(\brhop)-\bup^2))\\
\ =\ & B\bup +\bup(x)\partial_x(c^2-\bup^2) +\partial_x\bup(x)(c^2-\bup^2)\\
\ =\  & \bup\left[ \frac{a'}{a} p'(\brhop)-p''(\brhop)\partial_x\brhop
+\partial_x\left(\frac{M^2}{(\partial_x\bar{\psi}_+)^2}\right)\right]\\
&+\bup \partial_x
\left(p'(\brhop)-\frac{M^2}{(\partial_x\bar{\psi}_+)^2}\right) +
\partial_x\bup
(p'(\rho)-\bup^2)\\
\ =\
   & \bup \frac{a'}{a}p'(\brhop)+\frac{d
}{dx}\left(\frac{M}{a\brhop}\right)(p'(\brhop)-\bup^2).
\end{aligned}
\end{equation}
The momentum equation can be written as
\begin{equation*}
\partial_x\left(\frac{M^2}{a(x)\brhop}\right)+a(x)p'(\brhop)\brhop'=0,
\qquad
{\rm so,}
\qquad
\left(p'(\brhop)-\frac{M^2}{(a\brhop)^2}\right)\brhop'
=\frac{M^2}{a^3 \brhop} a'.
\end{equation*}
Therefore,
\begin{equation}\label{term2}
\bup \frac{a'}{a}p'(\brhop)+\frac{d }{dx}\left(\frac{M}{a\brhop}\right)(p'(\brhop)-\bup^2)=0.
\end{equation}

The boundary terms in (\ref{energy1}) are
\begin{equation}\label{bdyterm}
\begin{aligned}
\bup^2(x)(\partial_t   &  \Psi)^2    \Big|_{x=x_0}^{x=L}
\ +\
\bup (\bup^2-c^2)\partial_x\Psi\partial_t\Psi\Big|_{x=x_0}^{x=L}\\
&\ =\  \bup^2 (\partial_t\Psi)^2(\tau, L)+\bup^2(x_0)(\partial_t\Psi)^2(\tau, x_0)+\frac{a'}{a}\frac{\bup^2 \bum}{2 }\partial_t \Psi^2(\tau, x_0)\\
&\ =\  \bup^2 (\partial_t\Psi)^2(\tau, L)+\bup^2(x_0)(\partial_t\Psi)^2(\tau, x_0)+\partial_t\left(\left(\frac{a'\bup^2\bum \Psi^2}{2a}\right)(\tau, x_0)\right).
\end{aligned}
\end{equation}
The lemma follows from (\ref{energy1}), (\ref{term1}), (\ref{term2}), and (\ref{bdyterm}).
\end{pf}


The  estimate
 in Lemma \ref{keylem} and the
 method of  \cite{LRXX} imply  that solutions of the linear
 problem \eqref{lzeropb}
 decay exponentially to zero. This
 in turn  implies the exponential decay of
 solutions of the nonlinear problem.
 We describe the form that  these arguments
 take in the present context.

 \begin{lem}
 \label{lem:RT}
 If $a'$ satisfies (\ref{stabilitycondition})
 then there exist constants $\lambda_0>0$ and $C>0$
so that
solutions $\Psi$ of the problem (\ref{lzeropb})
satisfy,
\begin{equation*}
E(\Psi,t)
\ \leq\
 C\,e^{-\lambda_0 t}\, E (\Psi, 0)
\end{equation*}
\end{lem}

\begin{pf}
Step 1. Rauch-Taylor  estimates.
Thanks to the boundary condition at $x=x_0$,  the estimate \eqref{energyED} implies
\begin{equation*}
E(\Psi, t)+C_1\int_{0}^t (\Psi_t^2+\Psi_x^2)(s, x_0)\, ds
\ \leq\
 E(\Psi, 0)+C_2 \int_{0}^t \Psi^2 (s,x_0)\, ds\,.
\end{equation*}
The argument in \cite{Rauch} and the details in \cite[Appendix]{LRXX}, one can show that there exist
$ T>0$ and $\delta\in (0,  T/4)$ such that
\begin{equation} \label{estimateRT}
\begin{aligned}
\int_{0}^{T} (\Psi_t^2+\Psi_x^2)(t, x_0)\,dt
\ \geq\
 & \delta E(\Psi,  T)-C_3\int_0^{ T} \Psi^2(t, x_0)\,dt\,.
\end{aligned}
\end{equation}
Combining (\ref{energyED}) and (\ref{estimateRT}) yields,
\begin{equation}\label{energycpt}
(1+C_4)E(\Psi,  T)\ \leq\  E(\Psi, 0)+C_5\int_{0}^{ T} \Psi^2 (t, x_0)\, dt \, ,
\end{equation}
for some positive constants $C_4$ and $C_5$, independent of $t$.

Step 2. Spectrum of the evolution operator.
Define a new norm $\|\cdot\|_{\bX}$ for the function $h=(h_1,h_2) \in H^1\times L^2([x_0,L])$,
\begin{equation*}
\|h\|_{\bX}^2= \frac{a' \bup^2 \bum}{a}(x_0)|h_1|^2(x_0)+\int_{x_0}^L \bup \left\{|h_2|^2+
(p'(\bar{\rho}_+)-\bup^2)|h_1'|^2 \right\}(x) dx.
\end{equation*}
The associated complex Hilbert space will be denoted by $(\bX, \|\cdot\|_{\bX})$.
Define the solution operator $S_t: \bX\mapsto \bX$ as
\begin{equation*}
S_t(h)=(\Psi(t,\cdot), \Psi_t(t,\cdot))
\end{equation*}
where $\Psi$ is the solution of the problem (\ref{lzeropb}) with the initial data $h=(h_1,h_2)$.
Applying the Lemma on page 81 in \cite{Rauch2}, there are at most a finite set of generalized eigenvalues for the operator $S_{ T}$ in the annulus $\{\frac{1}{1+C_6}<|z|\leq 1\}\subset\mathbb C$, each with finite multiplicity.

Step 3. Refined estimate for the spectrum of $S_{ T}$. We show that the spectrum does
not touch the unit circle.
Otherwise there would exist $\omega\in \mathbb{R}$ and $V\in \bX$ such that
\begin{equation*}
(S_{ T}-e^{i\omega} I)V=0.
\end{equation*}
Note that the identity (\ref{energyED}) still holds in the complex setting if we replace the square terms in $E$ and $D$ by the square of modulus. Thus
\begin{equation*}
E(\Psi, 0)- E(\Psi, n T)=n D(\Psi,  T).
\end{equation*}
Since  $E(\Psi, n  T)$ and $E(\Psi, 0)$ are both positive and finite, it follows that
\begin{equation}\label{zerodsp}
D(\Psi,  T)=0.
\end{equation}
Therefore $E(\Psi, t)=E(\Psi, 0)$, for all $t$. Let $\mBV:=ker(S_{ T}-e^{i\omega}I)$.
Note that the coefficients in the problem (\ref{lzeropb}) do not depend on $t$, so
\[
(S_{ T}-e^{i\omega}I)S_t=S_t(S_{ T}-e^{i\omega}I)
\]
In particular, $S(t)\mBV\subset\mBV$ for any $t$, so $\mBV$ is invariant with respect to $S_t$.
Therefore
 $S_t|_{\mBV}$ is a semigroup on a finite dimensional subspace.
This yields that $S(t)|_{\mBV}=e^{tA}$
for some $A\in {\rm Hom}\,\mBV$.
The definition of $\mBV$ implies that $e^{ T}A =e^{i\omega t}I$  so the spectrum
of $A$ is purely imaginary.
Choose an eigenvector $w$ of $A$, $Aw=\lambda w$, $\lambda\in i\beta\in i\mBR$
Then
$S_t w=e^{\lambda t}w$. From
$S_t w=e^{i\beta t}w$ and  (\ref{zerodsp}) it follows that $w(x_0)=0$.
The boundary condition at $x_0$ then yields $w'(x_0)=0$.
The uniqueness of solutions of the linear homogeneous second order
ordinary differential equation satisfied by eigenfunctions implies that $w=0$.
This contradicts the assumption that $w$ is an eigenvector hence not equal to zero.

The contradiction shows that there exists $0<\beta_0<1$ so that
the spectrum,
$\mathfrak{\sigma}(S_{ T} )$, of $S_{ T}$ satisfies
$\mathfrak{\sigma}(S_T)
\, \subset\,
\{|z|\leq {\beta_0}\}
$.
The formula for the spectral radius implies that
$\|S(nT)\|=\|S(T)^n\|$ decays exponentially as $n\to \infty$.
This is equivalent to the assertion  Lemma \ref{lem:RT}.
\end{pf}

\begin{cor}
\label{cordecay}
Assume that $a'$ satisfies (\ref{stabilitycondition}) and $0\le k\in {\mathbb N}$.
Define
 \begin{equation*}
E_k(\Psi, t)\ :=\ \sum_{m=0}^{k}E(\partial_t^m\Psi, t)\,.
\end{equation*}
Then,
with $\lambda_0$ from Lemma \ref{lem:RT},
 solutions $\Psi$  of
 the linearized problem \eqref{lzeropb} satisfy
 \begin{equation*}
E_k(\Psi,t)
\ \leq\
 C\,e^{-\lambda_0 t}\, E_k (\Psi, 0)\,,
 \qquad
 {\rm and},
\end{equation*}
\begin{equation*}
\int_0^{\infty} e^{\lambda_0 t / 4}
\
\sum_{l=0}^{k}
\big(
 |\partial_t^l \Psi|^2(t, x_0)+|\partial_t^l \Psi|^2(t, L)
 \big)\ dt
\ \leq\
 C\, E_k(\Psi, 0).
\end{equation*}
\end{cor}

\begin{pf}
Since the coefficients of the equation and the boundary conditions
are independent of $t$, applying this estimate to the solution
$\partial_t^m\Psi$ yields
\begin{equation*}
E(\partial_t^m\Psi, t)
\ \leq\
C\, e^{-\lambda_0 t}\, E(\partial_t^m\Psi, 0).
\end{equation*}
Summing
on $m$ yield the first estimate of the Corollary.
The second follows from \eqref{energyED}, Sobolev imbedding, and the  first.
\end{pf}


\subsection{Uniform A Priori Estimates}\label{secuap}

The existence of local-in-time solutions for the problem \eqref{nonlinearpb}
is proved as   in \cite{LY}.
 In order to get global existence for
 the nonlinear  problem  \eqref{nonlinearpb},
 it suffces  to prove global {\it a priori}
 estimates for  \eqref{nonlinearpb} supplemented with initial data
 \begin{equation}\label{NLIC}
 \Psi(0, x)=h_1(x),\quad \Psi_t(0, x) =h_2(x).
 \end{equation}

With $C,\lambda_0$ from Lemma \ref{lem:RT}, choose $T>0$ so that
\begin{equation}
\label{eq:alphazero}
\alpha_0\ :=\ C\, e^{-\lambda_0 T}<1\,.
\end{equation}
For $t\geq T$  and  integer $k\ge 15$,
introduce
\begin{equation*}
\lbn (\Psi, \s)\rbn=\tn(Y, \s)\tn +\ttn(Y, \s)\ttn
\end{equation*}
where
\begin{equation*}
\begin{aligned}
\tn (\Psi, &  \sigma)\tn :=\sup_{\tau\in [0, t)}\sum_{0\leq m\leq k-6}\left(\sum_{0\leq l\leq m} e^{\frac{\lambda \tau}{16}}\|
\partial_t^l \partial_x^{m-l}\Psi(\tau, \cdot) \|_{L^{\infty}([x_0, L])} +e^{\frac{\lambda \tau}{16}}\left|\frac{d^m \s}{dt^m}\right|(\tau)\right)
\end{aligned}
\end{equation*}
and
\begin{equation*}
\begin{aligned}
&\ttn(\Psi, \s)\ttn =   \sup_{0\leq \tau \leq t} \left( \sum_{\substack{ 0\leq l \leq m,\\0\leq m\leq k  } }  \|\partial_t^l \partial_x^{m-l}\Psi(\tau, \cdot) \|_{L^{2}([x_0, L])}   +  \sum_{ 0\leq l \leq k}   \|\partial_t^l \partial_x^{k+1-l}\Psi(\tau, \cdot) \|_{L^{2}([x_0, L])}  \right)\\
&  + \sup_{0\leq \tau \leq t}   \|\partial_t^{k+1}\Psi(\tau, \cdot)-\frac{d^{k+1}\s}{d t^{k+1}}q_1(\cdot, \s)\partial_x \Psi(\tau, \cdot) \|_{L^{2}([x_0, L])} \\
& +\sum_{\substack{0\leq
l\leq m,\\ 0\leq m\leq k+1}} (\|\partial_t^l\partial_x^{m-l}
\Psi(\cdot,x_0)\|_{L^2[0, t]} +\|\partial_t^l\partial_x^{m-l}
\Psi(\cdot,L)\|_{L^2[0, t]})+\sum_{0\leq m \leq k+1}\left\|\frac{d^m\sigma}{dt^m}\right\|_{L^2[0, t]},
\end{aligned}
\end{equation*}
with a positive constant $\lambda$ to be defined in (\ref{deflam})  below.

In this subsection, when there are no specific indications to the contrary,  we
 assume $a_{ij}$,  $b_i$, and $g$ are functions of $(x, \Psi, \nabla \Psi, \sigma, \dot{\sigma})$, $d_1$ and $e_1$  are functions of $(\Psi, \Psi_t)$.
 Define
 \begin{equation}\label{defmcE}
\begin{aligned}
\mcE(\Phi, t)
\ :=\
\frac{e_1 \Phi^2}{d_1\brhop}(t,x_0)
\ +\
 \int_{x_0}^L \bup\{(\partial_t \Phi)^2 - a_{11}(\partial_x\Phi)^2\} (t,x)\ dx.
\end{aligned}
\end{equation}
Furthermore, for any $l\in \mathbb{N}$ and given $\Psi$ and $\s$ such that $\lbn (\Psi, \s)\rbn<\infty$, define
$$
\mcE_l(\Phi, t) :=
\sum_{m=0}^{l}\mcE(\partial_t^m\Phi, t), \quad  {\rm and},
\quad
\hat{\mcE}_l(\Phi, t)
\ :=\
\mcE_{l-1}(\Phi, t)+\mcE_0(\partial_t^l \Phi-q_1(x, \s) \Psi_x\frac{d^l \s}{dt^l}, t)\,.
$$
It is easy to see that  if $\lbn (\Psi, \sigma)\rbn \leq \epsilon$ for sufficiently small $\epsilon>0$, then
\begin{equation*}
\mcE(\Phi, t)(t)
\ \ge\  C\, \int_{x_0}^L (\partial_t \Phi)^2 + (\partial_x\Phi)^2 +\Phi^2)(t,x)\ dx
\end{equation*}
for some constant $C>0$ independent of $t$.

\begin{prop}\label{propapest}
Assume that $a'$ satisfies (\ref{stabilitycondition}). There
exists an $\epsilon_0>0$ so that for any $0<\epsilon <\epsilon_0$,
if  $(\Psi, \s)$ is a smooth solution of the problem (\ref{nonlinearpb}) and \eqref{NLIC}
satisfying
\begin{equation*}
|\s_0|
\ +\
 \|h_1\|_{H^{k+2}}
 \ +\
 \|h_2\|_{H^{k+1}}
 \ \leq\
 \epsilon^2
 \ \leq\
  \epsilon_0^2,
  \quad
  {\rm and},
  \quad
  \lbn (\Psi, \s)\rbn\ \le\  \epsilon,
\end{equation*}
 then
\begin{equation}\label{apest}
 \lbn (\Psi,\s)\rbn\  \leq\  \epsilon/2\, .
\end{equation}
\end{prop}
\begin{pf}
The proof has four steps.

Step 1.  Lower order energy estimate. Define
\begin{equation*}
\mD(\Phi, t)
\ :=\
\int_0^t -\bup (a_{11}d_1 +a_{01}) (\partial_t \Phi)^2(\tau, x_0) d\tau  + \int_0^t \bup a_{01}(\partial_t \Phi)^2(\tau, L) d\tau\,,
\end{equation*}
and
$\mD_m(\Phi, t) =\sum_{l=0}^m \mD(\partial_t^l\Phi, t)$.
Taking the $m$-th ($0\leq m \leq k-1$) order derivative
of the equation (\ref{nonlinearpb}) with respect to
$t$, then multiplying the both sides by $\bup \partial_t^{m+1}\Psi$ and
integrating on $\Omega=:[0, t]\times [x_0, L]$, noticing that $a_{00}=1$
yields
\begin{equation}\label{nllenergy}
\begin{aligned}
\mcE_m(\Psi, t)\  &+ \
 \mD_m(\Psi, t)
 \ \leq\
  \mcE_m(\Psi, 0) \\
&  +\   C\,\lbn(\Psi, \sigma)\rbn \left[ \int_0^t
 e^{-\frac{\lambda \tau}{64}}
 \left[ \mcE_m(\Psi, \tau) +   \sum_{l=0}^{m+2} \left|\frac{d^l \s}{dt^l}\right|^2  \right]  d\tau +\mD_m(\Psi, t)\right],
\end{aligned}
\end{equation}
for $m=0, 1, \cdots ,k-1$.

Step 2. The highest order energy estimates. Take the $k$-th order derivative for the equation (\ref{nonlinearpb}) with respect to
$t$ to find,
\begin{equation*}
\begin{aligned}
&\mL(x, \Psi, \sigma)\partial^k_t \Psi
\ =\ \mF_k(x, \Psi, \sigma)+ \frac{d^{k+2}\s}{dt^{k+2}}q_1(x, \s)\Psi_x +
\sum_{0\leq l \leq k-1}C_k^l\frac{d^{l+2}\s}{dt^{l+2}}\partial_t^{k-l}
\left(q_1(x, \s)\Psi_x\right),
\end{aligned}
\end{equation*}
where
\begin{equation*}
\mF_k(x, \Psi, \s) \ :=\  \sum_{1\leq l \leq k}C_k^l\left\{-\sum_{i,j=0}^1\partial_t^{l}a_{ij}\partial_{ij}\partial^{k-l}_tY- \sum_{i=0}^1\partial_t^{l}b_{i}\partial_{i}\partial^{k-l}_tY \right\}.
\end{equation*}
In order to handle the term $\frac{d^{k+2}\sigma}{dt^{k+2}}$,
study the equation for $\check{\Psi} =\partial_t^k \Psi-q_1(x, \s) \Psi_x \frac{d^k \s}{dt^k}$,
\begin{equation}\label{cptnlheq}
{\mL}(x, \Psi, \sigma)\check{\Psi}=\mF_k(x,\Psi, \sigma)+\check{\mF}(x,\Psi, \sigma),
\end{equation}
where
\begin{equation*}
\begin{aligned}
&\check{\mF}(x,\Psi, \sigma)
=\sum_{0\leq l \leq k-1}C_k^l\frac{d^{l+2}\s}{dt^{l+2}}\partial_t^{k-l}
(q_1(x, \s)\Psi_x) - \frac{d^{k+1}\s}{dt^{k+1}}
\partial_t\left(q_1(x, \s)\Psi_x\right) \\
&\quad  -\frac{d^{k}\s}{dt^{k}}
\partial_t^2\left(q_1(x, \s)\Psi_x\right)
-\left(2a_{01}\partial_t\partial_x+a_{11}\partial_x^2
+\sum_{i=0}^1b_i\partial_i \right) \left( \frac{d^{k}\s}{dt^{k}}
q_1(x, \s)\Psi_x\right).
\end{aligned}
\end{equation*}
Multiplying  both sides of (\ref{cptnlheq}) by $\bup\partial_t\check{\Psi}$ and
integrating on $\Omega$ yields
\begin{equation}\label{nlhenergy}
\begin{aligned}
\mcE(\check{\Psi}, t) \ +\ & \mD (\check{\Psi}, t)\ \leq\ \mcE(\check{\Psi}, 0)\\
&+C\lbn (\Psi, \sigma)\rbn \left[
\int_0^t e^{-\frac{\lambda \tau}{64}} \left[\hat\mcE_k(\Psi, \tau) +\sum_{l=0}^{k+1} \left|\frac{d^{l}\s}{dt^{l}}\right|^2 \right]d\tau +\hmD_k(\Psi, t)\right],
\end{aligned}
\end{equation}
where $\hmD_k$ is defined as follows
\begin{equation*}
\hat{\mD}_k(\Phi, t)=\mD_{k-1}(\Phi, t) +\mD_0(\partial_t^k \Phi-q_1(x, \s)\Psi_x\frac{d^k \s}{dt^k}, t).
\end{equation*}
Adding  the estimates (\ref{nllenergy}) and (\ref{nlhenergy}) yields
\begin{equation}\label{nlfullenergy}
\begin{aligned}
\hat{\mcE}_k(\Psi, t)\ +\ & \hat{\mD}_k(\Psi, t)\ \leq\
\hat{\mcE}_k(\Psi, 0)\\
&+\ C\,\lbn(\Psi, \sigma)\rbn \left\{
\int_0^t  e^{-\frac{\lambda \tau}{64}} \left[\hat\mcE_k(\Psi, \tau) +\sum_{l=0}^{k+1}\left|\frac{d^{l}\s}{dt^{l}}\right|^2 \right] d\tau +\hmD_k (\Psi, t)\right\}.
\end{aligned}
\end{equation}


Step 3. Boundedness  of the energy. Differentiating  the equation for the shock front
\begin{equation}\label{eqshockfront}
\s(t) =\mA_3(\Psi(t, x_0))
\end{equation}
with respect to $t$,  yields
\begin{equation*}\label{estsspeed}
\sum_{l=0}^{k+1} \left|\frac{d^l\s}{dt^l} \right|^2 (\tau)
\  \leq \
C\left( |\s(\tau)|^2 +\sum_{l=0}^{k+1}
|\partial_t^l \Psi(\tau, x_0)|^2   \right).
\end{equation*}
Using \eqref{eqshockfront} again yields
\begin{equation*}\label{estsig}
\sum_{l=0}^{k+1} \left|\frac{d^l\s}{dt^l} \right|^2\  \leq\  C \left(
\sum_{l=0}^{k} |\partial_t^l \Psi(\tau, x_0)|^2 +|\partial_t \check{\Psi}|^2 \right).
\end{equation*}
Thus
\begin{equation*}
\int_0^t e^{-\frac{\lambda \tau}{64}} \sum_{l=0}^{k+1}
\left|\frac{d^{l}\s}{dt^{l}}\right|^2 \, d\tau
\  \leq\
 C \left(\int_0^t
e^{-\frac{\lambda \tau}{64}}\mcE_0(\Psi, \tau; \Psi, \s) \,d\tau +\hat{\mD}_k(\Psi,t;
\Psi, \s)\right)
\end{equation*}
Therefore,  the energy estimate (\ref{nlfullenergy})  is equivalent to
\begin{equation*}
\begin{aligned}
&\hat{\mcE}_k(\Psi, t)+ \hat{\mD}_k(\Psi, t)
\ \leq\
 \hat{\mcE}_k(\Psi, 0)+C\lbn(\Psi, \sigma)\rbn \left ( \hat{\mD}_k(\Psi, \tau ) + \int_0^t  e^{-\frac{\lambda \tau}{64}}
\hat\mcE_k(\Psi, \tau)\, d\tau   \right).
\end{aligned}
\end{equation*}
If $\lbn (\Psi, \s)\rbn \leq  \epsilon$, then
\begin{equation*}
\hat{\mcE}_k(\Psi, t)\ +\   \hat{\mD}_k(\Psi, t)
\ \leq\
 2 \, \hat{\mcE}_k(\Psi, 0)\ \leq\
 C\,\epsilon^4.
\end{equation*}
This yields
\begin{equation}
\label{energybound}
\begin{aligned}
\ttn(\Psi, \s)\ttn
\ \leq \
 C\left(\sup_{0\leq \tau\leq t}\hmcE^{1/2}_k(\Psi, \tau) +\hmD_k^{1/2}(\Psi, t) \right)
\ \leq\
 C\epsilon^{2}
 \ \leq \
 \frac{\epsilon}{4}\,\,.
\end{aligned}
\end{equation}


Step  4. Decay of the  lower energy and  the  shock position.  The basic idea to get the
decay is to control the deviation of the solution $\Psi$ of the nonlinear problem  \eqref{nonlinearpb}
from that of the linearized problem (\ref{lzeropb}) (denoted by $\bPsi$). The contraction of the energy for $\bPsi$ will also yields the contraction of the energy for $\Psi$. This gives the decay of the lower energy of $\Psi$. The decay of the shock position is a consequence of the governing equation for the shock front and the decay of the lower energy.

At time $\tau=t_0$, we can choose $\bar{h}_1\in H^{k}$ and $\bar{h}_2\in H^{k-1}$ so
 that there exists a solution $\bPsi\in C^{k-1-i}([t_0,\infty); H^{i}([x_0, L]))$ of the linear problem (\ref{lzeropb}) satisfying $\bPsi(t_0, \cdot) =\bar{h}_1$ and $\bPsi_t(t_0, \cdot) =\bar{h}_2$, and $\bPsi$ satisfies
\begin{equation}\label{devIC1}
\sum_{l=0}^{k-1}\sum_{i=0}^l \|\partial_t^i \partial_x^{l-i}\bPsi(t_0, \cdot)\|_{L^2[x_0, L]}
\ \leq\
C\,\lbn (\Psi, \s)\rbn
\end{equation}
for some uniform constant $C$, and
\begin{equation}\label{devIC2}
\hmcE_{k-4}(\Psi-\bPsi, t_0)
\  \leq \ C \, \lbn (\Psi, \s)\rbn\  \hmcE_{k-4}(\Psi, t_0).
\end{equation}

As that in Steps 1 and 2, energy estimates for the equation for $\Psi-\bPsi$ gives
\begin{equation*}
\begin{aligned}
&\hmcE_{k-4}(\Psi-\bPsi, t_0+T) +\hmD_{k-4}(\Psi-\bPsi, t_0+T) - \hmD_{k-4}(\Psi-\bPsi, t_0) \\
\ \leq\
 & \hmcE_{k-4}(\Psi-\bPsi, t_0) +C \lbn (\Psi, \s)\rbn  \Bigg[ \int_{t_0}^{t_0+T}  \sum _{l=0}^{k-3} \left|\frac{d^{l}\s}{dt^l }\right|^2 + \hmcE_{k-4}(\Psi-\bPsi, \tau)\  d\tau \\
&\quad  +  \int_{t_0}^{t_0+T} \hmcE_{k-4}^{1/2}(\Psi, \tau) \hmcE_{k-4}^{1/2}(\Psi-\bPsi, \tau)   d\tau +( \hmD_{k-4}(\Psi, t_0+T) - \hmD_{k-4}(\Psi, t_0))\\
&\quad + (\hmD_{k-4}(\Psi-\bPsi, t_0+T) - \hmD_{k-4}(\Psi-\bPsi, t_0)) \Bigg].
\end{aligned}
\end{equation*}

Using the contraction of energy for $\bPsi$ and noting that the deviation of $\Psi$ and $\bPsi$ at $t_0$ is of higher order (cf. \eqref{devIC2}), one has
\begin{equation}\label{declowenergy}
\begin{aligned}
\frac{34+30\alpha_0}{64} \ \hmcE_{k-4}(\Psi, t_0+T)
\ \leq\
 \frac{2+30\alpha_0}{32}\  \hmcE_{k-4}(\Psi, t_0),
\end{aligned}
\end{equation}
if $\epsilon$ is sufficiently small.
As in  Corollary  \ref{cordecay},  it follows from (\ref{declowenergy}) that
\begin{equation*}
\hmcE_{k-4}(\Psi_t, t)+\s^2(t)
\ \leq\
 C\, \big(
 \tmcE_{k-4}(\Psi, 0)+\s^2(0)
 \big) \, e^{-2\lambda t},
\end{equation*}
where with $\alpha_0$ from \eqref{eq:alphazero},
\begin{equation}\label{deflam}
\lambda
\ :  =\
- \,
\frac{ \ln (1+\alpha)/2}{2T}, \qquad \text{and} \qquad
 \alpha \ :=\
 \frac{2+30\alpha_0}{17+15\alpha_0}.
\end{equation}
 Thus
\begin{equation}\label{maxnormdecay}
\sum_{l=0}^{k-6} \left[\left|\frac{d^l \s}{dt^l}\right| +\|\Psi(t,\cdot)\|_{L^{\infty}[x_0, L]} \right] \leq C\epsilon^2 e^{-\lambda t}.
\end{equation}
Combining (\ref{energybound}) and (\ref{maxnormdecay}), one has (\ref{apest}).
This finishes the proof of the Proposition \ref{propapest}.
\end{pf}

Once one has the Proposition \ref{propapest}, Theorem \ref{Thmdystability} follows from the standard continuation argument and local existence result.

\bigskip


{\bf Acknowledgments:} This work was initiated when Xie worked at University of Michigan and was completed when he visited The
Institute of Mathematical Sciences, The Chinese University of Hong
Kong. He gratefully acknowleges the both institutions. Rauch's
research is partially supported by  NSF grant DMS-0807600, and
Xin's research is partially supported by Hong Kong RGC Earmarked
Research Grants CUHK 4040/06P, CUHK 4042/08P, and a focus area scheme
grant at CUHK.



\begin{thebibliography}{99}



\bibitem{EmbidGM}
Pedro Embid, Jonathan Goodman,  and Andrew Majda,
{\it Multiple
steady states for {$1$}-{D} transonic flow}, SIAM J. Sci. Statist.
Comput.  {\bf 5} (1984), no. 1, 21--41.


\bibitem{Evans}
Lawrence C. Evans,  {\it Partial Differential Equations},  Graduate Studies in Mathematics, 19. American Mathematical Society, Providence, RI, 1998.



\bibitem{Ha}
Seung-Yeal Ha,  {\it $L\sp 1$ stability for systems of conservation laws with a nonresonant moving source},  SIAM J. Math. Anal. {\bf 33} (2001), no. 2, 411--439.

\bibitem{HaYang}
Seung-Yeal Ha and Tong Yang,  {\it $L\sp 1$ stability for systems of hyperbolic conservation laws with a resonant moving source}, SIAM J. Math. Anal.  {\bf 34} (2003), no. 5, 1226--1251


\bibitem{Kato}
Tosio Kato,  {\it Perturbation Theory for Linear Operators}, Reprint of the 1980 edition, Classics in Mathematics, Springer-Verlag, Berlin, 1995.

\bibitem{Lax}
Peter D. Lax,  {\it Functional Analysis}, Pure and Applied Mathematics (New York), Wiley-Interscience [John Wiley \& Sons], New York, 2002.





\bibitem{LY}
Ta Tsien Li and Wen Ci Yu, {\it Boundary value problems for quasilinear hyperbolic systems}, Duke University Mathematics Series, V. Duke University, Mathematics Department, Durham, NC, 1985.

\bibitem{Lien}
Wen-Ching Lien, {\it Hyperbolic conservation laws with a moving source}, Comm. Pure Appl. Math.  {\bf 52} (1999), no. 9, 1075--1098.

\bibitem{LiuARMA}
Tai-Ping Liu, {\it Transonic gas flow in a duct of varying area},
Arch. Rational Mech. Anal. {\bf 80} (1982), no. 1, 1--18.

\bibitem{LiuTP}
Tai-Ping Liu, {\it Nonlinear stability and instability of transonic
flows through a nozzle}, Comm. Math. Phys. {\bf 83} (1982), no. 2,
243--260.

\bibitem{LiuJMP}
Tai-Ping Liu,  {\it Nonlinear resonance for quasilinear hyperbolic equation}, J. Math. Phys.  {\bf 28} (1987), no. 11, 2593--2602.

\bibitem{LRXX}
Tao Luo, Jeffrey Rauch, Chunjing Xie, and Zhouping Xin, {\it
Stability of transonic shock solutions for Euler-Poisson equations}, to appear in
Arch. Rational Mech. Anal., arXiv:1008.0378,


\bibitem{Majda}
Andrew Majda,  {\it The existence of multidimensional shock fronts},  Mem. Amer. Math. Soc. {\bf 43} (1983), no. 281.



\bibitem{Metivier}
Guy M\'etivier, {\it Stability of multidimensional shocks}, Advances
in the theory of shock waves, 25--103, Progr. Nonlinear Differential
Equations Appl., 47, Birkh
user Boston, Boston, MA, 2001.

\bibitem{Rauch2}
Jeffrey Rauch,  {\it Qualitative behavior of dissipative wave equations on bounded domains},  Arch. Rational Mech. Anal.  {\bf 62} (1976),  no. 1, 77--85.



\bibitem{RauchM}
Jeffrey Rauch and Frank Massey, {\it Differentiability of solutions to hyperbolic initial-boundary value problems}, Trans. Amer. Math. Soc. {\bf 189} (1974), 303--318.

\bibitem{Rauch}
Jeffrey Rauch and Michael Taylor,  {\it Exponential decay of solutions to hyperbolic equations in bounded domains}, Indiana Univ. Math. J. {\bf 24} (1974), 79--86.





\bibitem{XinYinJDE}
 Zhouping Xin and Huicheng Yin, {\it The transonic shock in a nozzle, 2-D and 3-D complete Euler systems},
 J. Differential Equations {\bf 245} (2008), no. 4, 1014--1085.


\end{thebibliography}
\end{document}